\documentclass{article}
\usepackage{graphicx} 

\begin{document}

\title{Cayley's formula from middle school math}
\author{Victoria Feldman}
\maketitle

\begin{abstract}
The note contains a short elementary proof of Cayley's formula for labeled trees.
\end{abstract}

Cayley's formula states that there are $n^{n-2}$ trees on $n$ labeled vertices. This popular result was first proven by Borchardt with a determinantal technique \cite{Borchardt}. Cayley's proof  \cite{Cayley} employed the language of polynomials. The best-known derivation \cite{Prufer} by Pr\"ufer builds on a clever one-to-one correspondence. Several elegant proofs were collected in {\it Proofs from THE BOOK} \cite{THEBOOK}. Four proofs can be found in  {\it American Mathematical Monthly} \cite{AMM63,AMM16,AMM18,AMM24}. See Ref. \cite{AMM24} for a review of other derivations.

Many proofs are non-elementary. Others involve nontrivial combinatorics. We propose a proof with simple combinatorics. The main idea comes from using middle-school elementary algebra.

Recall that a tree is a graph with exactly one path without returns between any two vertices. Vertices are labeled with numbers $1$ to $n$. Two trees are identical if their edges connect the same pairs of labels. 

Let us build a tree with $n$ vertices. We will draw the tree and write the number of ways to label it. In the process, we will count trees of different types. Counting will be easy, but it will not be immediately obvious how to sum up all the cases. An elementary algebra trick will solve this difficulty.

We start by plotting a path from vertex $1$ to $2$. Let us assume there are $k_1\ge 2$ vertices on this path. If $k_1=2$, we write nothing. Otherwise we write the number of ways to label the $(k_1-2)$ intermediate vertices on the path: $(n-2)(n-3)\cdots (n-k_1+1)$.

After building the first path, we add a second path, which connects the first path to the vertex with the lowest label $x$ still unused. There are $k_1$ ways to choose from which vertex of the first path the second path will start.  Thus, we multiply the expression we wrote by $k_1$. Now that we have plotted a second path, let us assume that we have drawn a total of $k_2>k_1$ vertices. If $k_2>k_1+1$, we multiply the expression we wrote by the number of ways to label the second path: 
$(n-k_1-1)(n-k_1-2)\cdots(n-k_2+1)$. 

The third path we draw connects to the lowest label $y$ still unused. The total number of the vertices drawn is now $k_3$. The expression we wrote is multiplied by
$k_2(n-k_2-1)(n-k_2-2)\cdots (n-k_3+1)$.

We follow this pattern of adding edges and labeling vertices until we reach $n$ vertices. On each step, a new branch is added so that we plot $k_i>k_{i-1}$ vertices. 
The number we wrote becomes
$$\cdots (n-k_s+1)(k_s)(n-k_s-1)(n-k_s-2)\cdots(n-k_{s+1}+1)(k_{s+1})(n-k_{s+1}-1)$$
$$\times(n-k_{s+1}-2)\cdots (n-k_{p}+1)(k_p)(n-k_{p}-1)\cdots$$
There are $n-2$ factors in the expression.
To get the number of trees, these expressions should be summed over all choices of $2\le k_1<k_2<k_3<\dots k_{\rm final}=n$. 

The argument has been straightforward so far. It also seems unpromising since our sum looks unwieldy. Yet, the proof is almost done. Indeed, there are only two possibilities for the number $N_s$ at position $s$ from the left in the product. If there is $k_i=s+1$, then $N_s=s+1$. Otherwise $N_s=n-s-1$. Hence, by the distributive property of multiplication, the sum reduces to
$$
(2+n-2)(3+n-3)\ldots(n-1+n-[n-1])=n^{n-2}.
$$

To clarify this, let us consider the first few numbers of our long product.
So, at the first position,  $N_{1}$ either equals $2$ or $n-2$. If $N_1=2$ and we add all the products that start with $N_1=2$, we get 
$$2\cdot \sum N_2 \cdot N_3 \cdot\ldots \cdot N_{n-2}.$$
If $N_1=n-2$, the sum of the products is 
$$(n-2)\cdot \sum N_2 \cdot N_3 \ldots\cdot N_{n-2},$$
where the possible choices of $N_2, N_3,\ldots, N_{n-2}$  are the same in both sums since they only depend on $k_s>2$, and the same choices of the set of $k_s>2$ are allowed for
$k_1=2$ and $k_1>2$.
 Adding these two cases for $N_1$ gives us
$$(2)\cdot  \sum N_2 \cdot N_3 \cdot\ldots \cdot N_{n-2}
+(n-2)\cdot  \sum N_2 \cdot N_3 \cdot\ldots \cdot N_{n-2}=n  \sum N_2 \cdot N_3 \cdot\ldots \cdot N_{n-2}.$$
Each product in the sum $ \sum N_2 \cdot N_3 \cdot\ldots \cdot N_{n-2}$ starts with $N_2$. If $N_2=3$, we can write the sum of the products as 
$$3\cdot  \sum N_3 \cdot N_4 \cdot\ldots \cdot N_{n-2}.$$
If $N_2=n-3$, the sum of the products is
$$(n-3)\cdot \sum N_3 \cdot N_4 \cdot\ldots \cdot N_{n-2}.$$
After we add these two scenarios together, our total sum becomes
$$(n-2+2)(n-3+3) \sum N_3 \cdot N_4 \cdot\ldots \cdot N_{n-2}.$$
We repeat the same calculation for $N_3$ and so on.
The sum of the products collapses  and the number of trees becomes

$$(n-2+2)(n-3+3)\dots (1+n-1)=n^{n-2}.$$
\\

{\bf Acknowledgment} A useful discussion with R. E. Schwartz is gratefully acknowledged.


\begin{thebibliography}{99}
\bibitem{Borchardt} Borchardt, C. W. (1860). \"Uber eine Interpolationsformel f\"ur eine Art Symmetrischer Functionen und \"uber Deren Anwendung. Math. Abh. der Akademie der Wissenschaften zu Berlin, pp. 1–20. 

\bibitem{Cayley} Cayley, A. (1889). A theorem on trees. Q. J. Pure Appl. Math. 23: 376–378.

\bibitem{Prufer} Pr\"ufer, H. (1918). Neuer Beweis eines Satzes \"uber Permutationen. Arch. Math. Phys. 27: 742–744.

\bibitem{THEBOOK} Aigner, M., Ziegler, G. M. (2018). Proofs from THE BOOK, 6th ed. Berlin: Springer. doi.org/10.1007/978-3-662-57265-8.


\bibitem{AMM63} Moon, J. W. (1963).  Another Proof of Cayley's Formula for Counting Trees.
Amer. Math. Monthly. 70(8): 846-847.


\bibitem{AMM16} Avron, A., Dershowitz, N. (2016). Cayley’s formula: a page from the Book. Amer. Math. Monthly. 123(7): 699–700. 

\bibitem{AMM18} Shukla, A. B. (2018). A short proof of Cayley's tree formula. Amer. Math. Monthly. 125(1): 65-68.


\bibitem{AMM24} Zucker, M. (2024). From Abel’s Binomial Theorem to Cayley’s Tree Formula. Amer. Math. Monthly. 132(2): 165-169.
\end{thebibliography}
\end{document}